\documentclass{article}

\newcommand{\B}{{\cal B}}

\newcommand{\C}{{{\cal C}}}

\newcommand\one{\hbox{1\kern-2.4pt l }}

\newtheorem{theorem}{Theorem}[section]

\newtheorem{lemma}[theorem]{Lemma}
\newtheorem{corollary}[theorem]{Corollary}
\newtheorem{construction}[theorem]{Construction}
\def\whitebox{{\hbox{\hskip 1pt
    \vrule height 6pt depth 1.5pt
    \lower 1.5pt\vbox to 7.5pt{\hrule width
            3.2pt\vfill\hrule width 3.2pt}%
    \vrule height 6pt depth 1.5pt
    \hskip 1pt } }}
\def\qed{\ifhmode\allowbreak\else\nobreak\fi\hfill\quad\nobreak
         \whitebox\medbreak}

\begin{document}
\baselineskip 18pt

\title{On balanced $(Z_{4u}\times Z_{8v},\{4,5\},1)$ difference packings}{}

\author{{\small Hengming Zhao}\\
{\small  Department of Mathematics and Information Science}\\
{\small  Guangxi College of Education}\\
{\small  Nanning 530023, P. R. China}\\
{\small Rongcun Qin}\\
{\small  Xingjian College of Science and Liberal Art}\\
{\small  Guangxi University}\\
{\small  Nanning 530005, P. R. China}\\
{\small Dianhua Wu\thanks{Corresponding author.}}\\
{\small Department of Mathematics}\\
{\small Guangxi Normal University}\\
{\small Guilin 541004,   China}}

\date{}
\maketitle \noindent {\bf Abstract} Let $K$ be a set of positive integers and let $G$ be an additive
group. A $(G, K, 1)$ difference packing is a set of subsets of $G$
with sizes from $K$ whose list of differences covers every
element of $G$ at most once.
 It is balanced  if  the number of blocks of size $k\in K$ does not depend on $k$.
 In this paper, we determine a balanced $(Z_{4u}\times Z_{8v},{4,5},1)$
difference packing of the largest possible size whenever $uv$ is
odd. The corresponding optimal balanced  $(4u, 8v,\{4,5\},1)$ optical
 orthogonal signature pattern codes  are also obtained.
 

\section{Introduction}

\indent \ \ \ \  As usual, the {\em list of differences} of a subset $C$ of an additive group $G$ will be denoted by $\Delta C$.
  This is the multiset
of all differences $x-y$ with $(x, y)$  an ordered pair of distinct elements of $C$.
More generally, the list of differences of a set $\C$ of subsets of $G$ is the multiset
$\Delta \C=\bigcup\limits_{C\in \C}\Delta C$.

\indent \ \ \ \ Let  $G$ be an additive group and $K$  a set  of positive integers, a $(G, K, 1)$
difference packing (DP) is a set $\B$ of subsets of $G$ (blocks) with sizes from $K$
such that $\Delta \B$ (the list of differences of $\B$) covers every element of $G$ at most
once. The {\em difference leave} of a $(G, K, 1)$ difference packing  is the set of elements of $G$ not
covered by the list of differences of the packing.
The difference packing is balanced (BDP) if the number of blocks of size $k\in K$ does not depend on $k$.
 Thus, the number of blocks of a $(G,K,1)$-BDP is $b|K|$ for a suitable constant $b$ (of course the converse
 is not true in general).

\indent \ \ \ \  Let $I(G)$ be the set of all the involutions and the identity element of $G$.
It is  evident that $I(G)$ is always contained in the difference leave of a $(G,K,1)$-DP. Consequently,
if $b|K|$ is the number of blocks of a $(G,K,1)$-BDP, we have $b\sum\limits_{k \in K}(k^2-k) \leq |G \setminus I(G)|$.
It follows that
$$b\leq \biggl{\lfloor} {|G\setminus I(G)|\over\sum\limits_{k\in K}(k^2-k)} \biggl{\rfloor}.$$
Thus a  $(G, K, 1)$-BDP should be called {\em optimal } when, for every $k\in K$, the number
of its blocks of size $k$ is the right hand side of the above inequality. A $(G, k, 1)$-DP
is a  $(G, K, 1)$-BDP with $K=\{k\}$. An $(u\times v, k,1)$-DP is a $(G, k, 1)$-DP with $G=Z_u\times Z_v$.
Similarly, a $(Z_u\times Z_v, K, 1)$-BDP is denoted by  $(u\times v, K,1)$-BDP.

\indent \ \ \ \  Let $s|u$, and $t|v$.
We speak of an $(s\times t)$-regular $(u\times v,K,1)$-BDP to mean a $(G,K,1)$-BDP with
$G=Z_u\times Z_v$ and difference leave $H=(u/s)Z_u\times (v/t)Z_v$. So it is actually a balanced
$(G,H,K,1)$ difference family (see \cite{Bur2, ZQ}). In the following such a BDP will be denoted by
$(u\times v, s\times t, K,1)$-BDP.

\indent \ \ \ \  If $b|K|$ is the number of blocks of an $(u\times v, s\times t, K, 1)$-BDP, then we have
$$b=\frac{uv-st}{\sum\limits_{k\in K}(k^2-k)}.$$
Since $I(Z_u\times Z_v)\subseteq (u/s)Z_u\times (v/t)Z_v$,
then we have $1\leq |I(Z_u\times Z_v)|\leq st$. If $1\leq st\leq
\sum\limits_{k\in K}(k^2-k)$, then
 $$\left\lfloor {uv-|I(Z_u\times Z_v)|\over\sum\limits_{k\in K}(k^2-k)}\right\rfloor=
\left\lfloor {uv-st+(st-|I(Z_u\times
Z_v)|)\over\sum\limits_{k\in K}(k^2-k)}\right\rfloor=
{uv-st\over\sum\limits_{k\in K}(k^2-k)}.$$ Thus, an $(u\times v,
s\times t, K,1)$-BDP is optimal when $1\leq st\leq \sum\limits_{k\in
K}(k^2-k)$.

\indent \ \ \ \  The construction  of an optimal
$(G,K, 1)$-BDP becomes easier when the order of $G$ is divisible by
$\sum\limits_{k\in K}(k^2-k)/2$.
An optimal
$(m\times n,K,1)$-BDP is equivalent to an optimal
$(mn,K,1)$-BDP when $m$ and $n$ are coprime.
In this case, one can obtain many optimal $(m\times n,K,1)$-BDPs
 from the results on optimal  $(mn,K,1)$-BDPs with $|K|>1$ (see \cite{BWWFC}, \cite{DVR},
\cite{GW}, \cite{JWL}, \cite{WZF}, \cite{Yang},
\cite{Zhao}, \cite{ZWF} for some of the examples. The notion $(mn,K,1)$-CPs is used in some of these references).
However, when $m$ and $n$ are not coprime, the construction of optimal $(m\times n,K,1)$-BDPs
is difficult.  Optimal  $(3u\times 9v, \{3,4\},1)$-BDPs  were constructed for any pair
of positive integers $(u,v)$  with {\rm gcd}$(u,3)=1$ \cite{ZQ}. Optimal  $(2u\times 16v, \{4,5\},1)$-BDPs
 were constructed for $uv$ odd  \cite{LZQW}. To the authors' knowledge,
 no other work had been done on optimal $(m\times n, \{4,5\},1)$-BDPs when $m$ and $n$ are not coprime.
In this paper, optimal $(4u\times 8v, \{4,5\},1)$-BDPs will be  constructed for $uv$ odd.

\indent \ \ \ \  The concept of a strong difference family was introduced in \cite{B-SDF} to provide constructions
for relative difference families (see also \cite{Bura-SG}, \cite{B-Z}, \cite{M}). Some work also had been
 done on strong difference families recently \cite{CCFW}, \cite{CFW}, \cite{CFW1}, \cite{WFW}.
 In \cite{BBGT}, a strong difference family was  used to construct $3$-pyramidal KTSs (Kirkman triple systems).
 In \cite{B-H},  \cite{BYW},  a strong difference family was also used to construct partitioned difference families.

\indent \ \ \ \  Let $G$ be an additive group, $K$ a set of positive integers, $\lambda$ a positive integer.
A $(G, K, \lambda)$ strong difference family ($(G, K, \lambda)$-SDF for short) $\B$  is a
 family of multisets of $G$ (blocks)  with sizes from $K$ such that $\Delta \B$ (the list
 of differences of $\B$) covers every element of $G$  exactly $\lambda$ times.

\indent \ \ \ \   Let $a, \ b$ be positive integers, $p$ an odd prime.
One can construct a $(Z_a\times Z_b\times Z_p, Z_a\times Z_b\times \{0\},K,1)$-BDP  by
 using a suitable $(Z_a\times Z_b,K,\mu)$-SDF.
In this paper we give direct and explicit constructions
for $(Z_a\times Z_b\times Z_p, Z_a\times Z_b\times \{0\},\{4,5\},1)$-BDPs.
This will be realized with the implicit use of suitable
$(Z_a\times Z_b,\{4,5\},2)$- or $(Z_a\times Z_b,\{4,5\},4)$-SDFs according to whether
$p\equiv3$ or 1 (mod 4), respectively. We will give the constructions in Section 2.

\indent \ \ \ \ In this paper, by using direct and recursive
constructions, the following result is obtained.

\begin{theorem}\label{main1}
If $u>1$, $v>1$ are odd integers, then there exists an optimal
$(4u\times 8v,\{4,5\},1)$-BDP.
\end{theorem}

\section{Direct Constructions}

 \indent \ \ \ \ In this section, we give  direct constructions for
$(4\times gp,4\times g,\{4,5\},1)$-BDPs with $g=8,24$ and prime
$p>5$. Quadratic residues  will be used in the constructions.

\indent \ \ \ \ Let $p$ be an odd prime. In the sequel, we will always
assume that $\omega$ is a primitive element of $Z_p$. Then,
$C_0^2=\{\omega^{2i}|0\leq i\leq (p-3)/2\}$ is called the set of
quadratic residues (squares) of $Z_p$, and $C_1^2=\{\omega^{2i+1}|0\leq i\leq
(p-3)/2\}$ is the set of quadratic non-residues (non-squares) of $Z_p$. Let $\xi$
be the first   quadratic non-residue of $Z_p$, it is easy to see
that $\xi$ is a prime. In the following lemma, the first four
results are easy to be obtained by the law of quadratic reciprocity
in number theory, the last result is stated in \cite{BWWFC}.

\begin{lemma}
\label{le33} \  For a fixed odd prime $p$, we
have:\\
\noindent ${\rm (1)}$ \ $-1\in C_0^2$ if and only if $p\equiv1 \pmod
4$;\\
\noindent ${\rm (2)}$ \ $2\in C_0^2$ if and only if $p\equiv \pm1
\pmod 8$;\\
\noindent ${\rm (3)}$ \ $3\in C_0^2$ if and only if $p\equiv\pm1
\pmod{12}$;\\
\noindent ${\rm (4)}$ \ $5\in C_0^2$ if and only if $p\equiv\pm1
\pmod{10}$;\\
\noindent ${\rm (5)}$ If $p\equiv\pm1 \pmod 8$, then
$\{\xi-2,\xi-1,\xi+1\}\subset C_0^2$.
\end{lemma}

\indent \ \ \ \ Let $y$ and $y+1$ be two consecutive non-squares in
$Z_p$, and  let $z$ be the smallest integer in
$\{y+2,$ $y+3,\ldots,p-1\}$ such that $z$ is a square of
$Z_p$. Take $\theta=z-1$,  the following result is from \cite{Bur1}.

\begin{lemma}
\label{le33-1} \ If $p\geq5$ is a prime, then there is an element
$\theta\in C_1^2$  such that $\theta-1\in C_1^2$ and $\theta+1\in
C_0^2$
\end{lemma}

\indent \ \ \ \ In the sequel, we will always assume that $\theta$
satisfies the properties in Lemma \ref{le33-1}. Two constructions
are stated below.

{\bf C1:} $p\equiv3$ (mod 4), $p \not| \ ab$.

\indent \ \ \ \  Let $\sigma=\{S_1, ..., S_h\}$ be a $(Z_a\times Z_b, \{4,5\},2)$-SDF
and let ${\cal B}=\{B_1,...,B_h\}$ be a set of subsets of $Z_a\times Z_b\times Z_p$
such that the projection of $B _k$ on $Z_a\times Z_b$ is $S_k$ for $k=1,...,h$.
Considering that $\sigma$ is a $(Z_a\times Z_b, \{4,5\},2)$-SDF, it is clear that
$\Delta {\cal B}$ has the form
$\Delta {\cal B}=\bigcup\limits_{z\in Z_a\times Z_b} z\times L_{z}$
where each $L_z$ is a pair of elements of $Z_p\setminus\{0\}$,
say $L_z=\{\ell_z,\ell'_z\}$.
In the case that $\ell_z\ell'_z$ is a non-square for each $z\in Z_a\times Z_b$
we have $L_z\cdot C_0^2=Z_p\setminus\{0\}$ for each $z$ and hence
${\cal F}=\{(1,1,s)\cdot B_{k} \ | \ s\in C_0^2; 1\leq k\leq h\}$ is a
$(Z_a\times Z_b\times Z_p, Z_a\times Z_b\times \{0\}, \{4,5\},1)$-BDP.

\smallskip
{\bf C2:} $p\equiv1$ (mod 4), $p \not| \ ab$.

\indent \ \ \ \  Let $\sigma=\{S_1, ..., S_{2h}\}$ be a $(Z_a\times Z_b, \{4,5\},4)$-SDF.
Then take a set ${\cal B}=\{B_1, ..., B_{2h}\}$ of subsets of $Z_a\times Z_b\times Z_p$
with the projections of $B_k$ on $Z_a\times Z_b$ is $S_k$ for $1\leq k\leq 2h$.
Considering that $\sigma$ is a $(Z_a\times Z_b, \{4,5\},4)$-SDF, we have
$\Delta {\cal B}=\bigcup\limits_{z\in Z_a\times Z_b} z\times L_z$
with $L_z$ a quadruple of elements of $Z_p\setminus\{0\}$ for each $z\in Z_a\times Z_b$.
Assume that the $B_k$'s can be taken in such a way that each $L_z$ is of the
form $L_z=\{1,-1\}\cdot\{\ell_z,\ell'_z\}$ with $\ell_z\ell'_z\in \ $$C_1^2$.
In this case we have $L_z\cdot C_0^2/\{1, -1\}=Z_p\setminus\{0\}$ for each $z$ and hence
${\cal F}=\{(1,1,s)\cdot B_k \ | \ s\in C_0^2/\{1, -1\}; 1\leq k\leq 2h\}$ is a
$(Z_a\times Z_b\times Z_p, Z_a\times Z_b\times \{0\}, \{4,5\},1)$-BDP.

\indent \ \ \ \  In the proof of the next lemma, we will apply {\bf C1}
using

 $\sigma=\{\{(0,0),(0,5),(2,1), (3,1)\},\{(0,0),(0,3),(1,6),(2,0)\},\\
\ \indent \ \ \ \ \ \ \  \{(0,0),(0,4),(0,6),(1,3),(3,2)\},\{(0,0),(0,0),(0,1),(1,1),(2,6)\}\}$

as a $(Z_4\times Z_8, \{4,5\},2)$-SDF
for the case of  $p\equiv3$ (mod 4).
Instead, for the case of  $p\equiv1$ (mod 4), we essentially apply {\bf C2}
using a $(Z_4\times Z_8, \{4,5\},4)$-SDF  below.

$\sigma=\{\{(0,0),(0,0),(0,1),(0,1)\},\{(0,0),(2,0),(0,4),(2,4)\},\\
\ \indent \ \ \ \ \ \ \ \{(0,0),(0,3),(1,4),(3,5)\},\{(0,0),(0,6),(1,0),(3,3),(3,6)\}\},\\
\ \indent \ \ \ \ \ \ \  \{(0,0),(0,6),(1,5),(2,1),(3,7)\},\{(0,0),(0,3),(1,4),(3,5)\},\\
\ \indent \ \ \ \ \ \ \  \{(0,0),(0,6),(1,0),(3,3),(3,6)\},\{(0,0),(0,6),(1,5),(2,1),(3,7)\}\}.$

\begin{lemma}
\label{lm25}  If $p\geq 5$ is a prime, then there exists a $(4\times8p,4\times8,\{4,5\},1)$-BDP.
\end{lemma}
{\bf Proof}  \ Since $p\geq 5$ is a prime, we identify $Z_4\times
Z_{8p}$ with $Z_4\times Z_{8}\times Z_p$. The problem is split into
two cases depending on the values of $p$ modulo 4.

{\it 1st case $p\equiv3\pmod{4}$}.

\noindent Consider four subsets of $Z_{4}\times Z_{8}\times Z_p$ of
the following form.

\indent \ \ \ \ $B_1=\{(0,0,0),(0,5,\alpha_1),(2,1,\alpha_2),(3,1,\alpha_3)\}$, \\
\indent \ \ \ \ $B_2=\{(0,0,0),(0,3,\beta_1),(1,6,\beta_2),(2,0,\beta_3)\}$,\\
\indent \ \ \ \ $B_3=\{(0,0,0),(0,4,\gamma_1),(0,6,\gamma_2),(1,3,\gamma_3),(3,2,\gamma_4)\}$,\\
\indent \ \ \ \
$B_4=\{(0,0,0),(0,0,\delta_1),(0,1,\delta_2),(1,1,\delta_3),(2,6,\delta_4)\}$.

\noindent We have $\bigcup\limits_{i=1}^4\Delta
B_i=\bigcup\limits_{(i,j)\in Z_4\times
Z_{8}}\{(i,j)\}\times L_{(i,j)}$ where:

\indent \ \ \ \ $L_{(0,0)}=\{\delta_1,-\delta_1\}$;
$L_{(0,1)}=\{\delta_2,\delta_2-\delta_1\}$;
$L_{(0,2)}=\{\gamma_2-\gamma_1,-\gamma_2\}$;

\indent \ \ \ \ $L_{(0,3)}=\{\beta_1,-\alpha_1\}$;
$L_{(0,4)}=\{\gamma_1,-\gamma_1\}$;
$L_{(1,0)}=\{\delta_3-\delta_2,\alpha_3-\alpha_2\}$;

\indent \ \ \ \ $L_{(1,1)}=\{\delta_3-\delta_1,\delta_3\}$;
$L_{(1,2)}=\{\beta_3-\beta_2,\gamma_1-\gamma_4\}$;
$L_{(1,3)}=\{\gamma_3,\beta_2-\beta_1\}$;

\indent \ \ \ \
$L_{(1,4)}=\{\alpha_1-\alpha_3,\gamma_2-\gamma_4\}$;
$L_{(1,5)}=\{\delta_4-\delta_3,\gamma_3-\gamma_2\}$;
$L_{(1,6)}=\{\beta_2,-\gamma_4\}$;

\indent \ \ \ \ $L_{(1,7)}=\{\gamma_3-\gamma_1,-\alpha_3\}$;
$L_{(2,0)}=\{-\beta_3,\beta_3\}$;
$L_{(2,1)}=\{\gamma_3-\gamma_4,\alpha_2\}$;

\indent \ \ \ \ $L_{(2,2)}=\{-\delta_4,\delta_1-\delta_4\}$;
$L_{(2,3)}=\{\beta_1-\beta_3,\delta_2-\delta_4\}$;
$L_{(2,4)}=\{\alpha_2-\alpha_1, \alpha_1-\alpha_2\}$;

\indent \ \ \ \ $L_{(0,j)}=-L_{(0,8-j)}$,
$L_{(2,j)}=-L_{(2,8-j)}$, $5\leq j\leq7$;
$L_{(3,j)}=-L_{(1,8-j)}$, $0\leq j\leq7$.

\indent \ \ \ \ Let
$(\alpha_1,\alpha_2,\alpha_3,\beta_1,\beta_2,\beta_3,\gamma_1,\gamma_2,\gamma_3,\gamma_4,\delta_1,\delta_2,\delta_3,\delta_4)$
 be  as follows:

\indent \ \ \ \ $(3,-2,1,-1,3,1,1,2,3,-1,2,1,3,-1), \ \ \ \ \ \ \ \ {\rm for}  \ p\equiv7\pmod{24};$\\
\indent \ \ \ \  $(-1,1,3,-1,3,2,1,4,2,3,1,2,3,-1),   \ \ \ \ \ \ \ \ \  \ {\rm for} \   p\equiv11\pmod{24};$\\
\indent \ \ \ \ $(1,4,2,1,-3,3,3,2,6,4,2,1,3,-1),  \  \ \ \ \ \ \ \ \ \ \ \ \ {\rm for}  \  p\equiv19\pmod{24};$\\
\indent \ \ \ \ $(3,2,1,3,1,2,-1,2,1,3,4,1,2,3),   \ \ \ \ \ \ \ \ \ \ \ \ \ \ \   {\rm for} \   p\equiv23\pmod{24}.$

\indent \ \ \ \ Using Lemma \ref{le33}, one can readily check that
$ L_{(i,j)}$ has a square and a non-square  for $(i,j)\in
Z_4\times Z_{8}$. Set $\B=\{(1,1,c)\cdot B_i| i=1,2,3,4,c\in
C_0^2\}$, we have $\Delta\B=\bigcup\limits_{(i,j)\in Z_4\times
Z_{8}}\{(i,j)\}\times( L_{(i,j)}\cdot
C_0^2)=\bigcup\limits_{(i,j)\in Z_4\times
Z_{8}}\{(i,j)\}\times(Z_p\setminus\{0\})$. Thus $\B$ forms a
  $(4\times8p,4\times8,\{4,5\},1)$-BDP.

{\it 2nd case $p\equiv1\pmod{4}$}.

\indent \ \ \ \ When $p=5$, the blocks of a
$(4\times40,4\times8,\{4,5\},1)$-BDP is displayed below:

\{(0,0),(1,29),(1,37),(2,6),(2,28)\}, \{(0,0),(1,4),(2,17),(3,18)\},

\{(0,0),(0,16),(1,27),(3,24)\}, \{(0,0),(0,7),(0,21),(2,8)\},

 \{(0,0),(0,1),(0,28),(0,37),(3,34)\}, \{(0,0),(1,18),(2,2),(3,14)\},

\{(0,0),(0,2),(1,23),(2,11),(3,4)\}, \{(0,0),(0,23),(0,29),(2,7),(3,21)\}.

\indent \ \ \ \ When $p>5$, by Lemma \ref{le33-1}, there exists a non-square
$\theta$ of $Z_p$ such that $\theta-1\in C_1^2$, $\theta+1\in
C_0^2$. Consider eight subsets of $Z_{4}\times Z_{8}\times Z_p$ of
the following form.

\indent \ \ \ \ $B_1=\{(0,0,1),(0,0,-1),(0,1,\theta),(0,1,-\theta)\}$, \\
\indent \ \ \ \ $B_2=\{(0,0,0),(2,0,\alpha_1),(0,4,\alpha_2),(2,4,\alpha_3)\}$, \\
\indent \ \ \ \ $B_3=\{(0,0,0),(0,3,\beta_1),(1,4,\beta_2),(3,5,\beta_3)\}$,\\
\indent \ \ \ \ $B_4=\{(0,0,0),(0,6,\gamma_1),(1,0,\gamma_2),(3,3,\gamma_3),(3,6,\gamma_4)\}$,\\
\indent \ \ \ \
$B_5=\{(0,0,0),(0,6,\delta_1),(1,5,\delta_2),(2,1,\delta_3),(3,7,\delta_4)\}$,\\
\indent \ \ \ \ $B_{5+i}=(1,1,-1)\cdot B_{2+i}$ for $i=1,2,3$.

\indent \ \ \ \ Let
$(\alpha_1,\alpha_2,\alpha_3,\beta_1,\beta_2,\beta_3,\gamma_1,\gamma_2,\gamma_3,\gamma_4,\delta_1,\delta_2,\delta_3,\delta_4)$
 be  as follows:

\indent \ \ \ \
$(\xi+\xi^2,\xi^2,2\xi^2,\xi^2,\xi,-\xi,1,\xi,\xi^2,-\xi,\xi,-\xi,-1,1)$,
\  for $p\equiv1\pmod{8}$;\\
\indent \ \ \ \   $(1,6,3,-4,6,-2,1,3,2,5,2,1,-1,5)$,\ \ \ \ \ \ \ \ \ \ \ \ \ \ \     for \  $p\equiv29,101\pmod{120}$;\\
\indent \ \ \ \ $(1,2,5,1,5,-5,2,5,3,1,1,-1,5,2)$,\ \ \ \ \ \ \ \ \ \ \ \ \  for \   $p\equiv53,77\pmod{120}$;\\
\indent \ \ \ \ $(1,3,9,1,3,2,1,-2,4,2,2,-2,6,4)$,\ \ \ \ \ \ \ \ \ \ \ \ \  for \   $p\equiv13\pmod{24}$.

\indent \ \ \ \  If $p\equiv1\pmod8$, then 
$\{\xi-2,\xi-1,\xi+1\}\subset C_0^2$ by Lemma \ref{le33}. Let
$\B=\{(1,1,c)\cdot B_i|$ $ 1\leq i\leq8, $ $c\in C_0^2/\{-1,1\}\}$, it is
checked that $\B$ forms a
$(4\times8p,4\times8,\{4,5\},1)$-BDP.\qed

\begin{lemma}
\label{lm26}  If $p\equiv3\pmod4$ is a prime and $p\geq7$, then
there exists a  $(4\times24p,4\times24,\{4,5\},1)$-BDP.
\end{lemma}
{\bf Proof} \ Since $p\geq7$ is a prime, then {\rm gcd}$(24,p)=1$,
 $Z_4\times Z_{24p}$ is isomorphic to $Z_4\times Z_{24}\times Z_p$.
Consider twelve subsets of $Z_4\times Z_{24}\times Z_p$ of the
following form:

\indent \ \ \ \ $B_1=\{(0,0,0),(0,0,\alpha_1),(0,1,\alpha_2),(2,7,\alpha_3)\}$,\\
\indent \ \ \ \ $B_2=\{(0,0,0),(0,6,\beta_1),(0,12,\beta_2),(2,0,\beta_3)\}$,\\
\indent \ \ \ \ $B_3=\{(0,0,0),(1,3,1),(1,11,2),(1,22,3),(3,14,4)\}$, \\
\indent \ \ \ \ $B_4=\{(0,0,0),(2,4,1),(2,14,2),(3,5,3),(3,8,4)\}$, \\
\indent \ \ \ \ $B_5=\{(0,0,0),(0,15,1),(1,14,2),(1,21,3),(3,19,4)\}$,\\
\indent \ \ \ \ $B_6=\{(0,0,0),(0,2,1),(1,2,2),(3,17,3)\}$,\\
\indent \ \ \ \ $B_7=\{(0,0,0),(0,4,1),(1,12,2),(3,11,3)\}$,\\
\indent \ \ \ \ $B_{5+i}=(1,1,-1)\cdot B_i$, $i=3,4,5,6,7$.

\indent \ \ \ \ Let
$(\alpha_1,\alpha_2,\alpha_3,\beta_1,\beta_2,\beta_3)$  be as
follows:

\indent \ \ \ \ $(2,1,3,1,4,2), \ \ \ \ \ \ \   {\rm for} \ \ p\equiv7\pmod{24};$\\
\indent \ \ \ \  $(2,1,4,1,3,2),  \ \ \ \ \ \ \   {\rm for} \ p\equiv11\pmod{24};$\\
\indent \ \ \ \ $(2,1,3,1,3,4),   \ \ \ \ \ \ \  {\rm for} \   p\equiv19\pmod{24};$\\
\indent \ \ \ \ $(3,1,2,2,1,3),   \ \ \ \ \ \ \  {\rm for} \  p\equiv23\pmod{24}.$

\indent \ \ \ \ Set $\B=\{(1,1,c)\cdot B_i| 1\leq i\leq12,c\in
C_0^2\}$. Using Lemma \ref{le33}, it is easy to check that $\B$
forms a  $(4\times24p,4\times24,\{4,5\},1)$-BDP. \qed

\begin{lemma}
\label{lm27}  If $p\equiv1\pmod4$ is a prime and $p>5$, then there
exists a  $(4\times24p,4\times24,\{4,5\},$ $1)$-BDP.
\end{lemma}
{\bf Proof} \ We identify $Z_4\times Z_{24p}$ with $Z_4\times
Z_{24}\times Z_p$. By Lemma \ref{le33-1}, there exists a non-square
$\theta$ such that $\theta-1\in C_1^2$, $\theta+1\in C_0^2$.
Consider eleven  subsets of $Z_4\times Z_{24}\times Z_p$ of the
following form:

\indent \ \ \ \
$B_1=\{(0,0,1),(0,0,-1),(0,2,\theta),(0,2,-\theta)\}$,\\
\indent \ \ \ \
$B_2=\{(0,0,0),(0,12,\alpha_1),(2,0,\alpha_2),(2,12,\alpha_3)\}$,\\
\indent \ \ \ \ $B_3=\{(0,0,0),(0,1,1),(1,1,2),(1,5,3),(1,11,4)\}$,\\
\indent \ \ \ \
$B_4=\{(0,0,0),(0,9,1),(1,17,2),(2,23,3),(3,15,4)\}$,\\
\indent \ \ \ \
$B_{5}=\{(0,0,0),(0,5,\beta_1),(0,21,\beta_2),(1,19,\beta_3),(2,16,\beta_4)\}$,\\
\indent \ \ \ \
$B_{6}=\{(0,0,0),(1,2,\gamma_1),(1,15,\gamma_2),(3,11,\gamma_3),(3,22,\gamma_4)\}$,\\
\indent \ \ \ \
$B_{7}=\{(0,0,0),(0,3,\delta_1),(1,22,\delta_2),(2,21,\delta_3)\}$,\\
\indent \ \ \ \
$B_{8}=\{(0,0,0),(0,7,\epsilon_1),(1,21,\epsilon_2),(3,4,\epsilon_3)\}$,\\
\indent \ \ \ \
$B_{9}=\{(0,0,0),(0,5,\varepsilon_1),(2,18,\varepsilon_2),(3,17,\varepsilon_3)\}$,\\
\indent \ \ \ \
$B_{10}=\{(0,0,0),(1,12,\zeta_1),(3,9,\zeta_2),(3,17,\zeta_3)\}$,\\
\indent \ \ \ \
$B_{11}=\{(0,0,0),(0,7,\eta_1),(1,3,\eta_2),(2,16,\eta_3)\}$.

\indent \ \ \ \  Let
$(\alpha_1, \alpha_2, \alpha_3,\beta_1, \beta_2, \beta_3, \beta_4,\gamma_1, \gamma_2, \gamma_3, \gamma_4,\delta_1,\delta_2,\delta_3)$
 be  as follows:

\indent \ \ \ \
$(\xi^2,\xi^2+\xi,2\xi^2,1,-1,\xi,-\xi,\xi^2,\xi^3,-\xi^2,\xi,\xi,\xi^2+\xi,\xi^2)$,
 \ \ \ for $p\equiv1\pmod{8}$;\\
\indent \ \ \ \ $(1,2,10,1,2,3,4,1,4,6,2,1,2,3)$, \ \ \ \ \ \ \ \ \ \ \ \ \ \ \ \ \  for $p\equiv29,101\pmod{120}$;\\
\indent \ \ \ \
$(1,2,5,1,2,3,4,1,4,6,2,1,2,3)$,  \ \ \ \ \ \ \ \ \ \ \ \ \ \ \ \ \ \ \    for
$p\equiv53,77\pmod{120}$;\\
\indent \ \ \ \
$(1,4,6,1,2,3,4,1,3,5,2,-4,2,-2)$, \ \ \ \ \ \ \ \ \ \ \ \ \ \   for
$p\equiv13,37\pmod{120}$;\\
\indent \ \ \ \
$(1,3,9,1,2,3,5,1,3,5,2,-4,2,-2)$, \ \ \ \ \ \ \ \ \ \ \ \ \ \ \  for
$p\equiv61,109\pmod{120}$.

\indent \ \ \ \ Let
$(\epsilon_1, \epsilon_2, \epsilon_3,\varepsilon_1, \varepsilon_2,\varepsilon_3,\zeta_1, \zeta_2,\zeta_3,\eta_1, \eta_2,\eta_3)$
 be  as follows:

 \indent \ \ \ \
$(\xi+1,1,\xi,\xi,\xi^2,1,-\xi,\xi^2,\xi,\xi,\xi^2+\xi,\xi^2)$, \ \ \ \ \ \ \ for $p\equiv1\pmod{8}$;\\
\indent \ \ \ \ $(1,2,6,2,1,3,2,3,1,2,10,12 )$, \ \ \ \ \ \ \ \ \ \ \ \ \ \ \ \ \ \ \ \ \  for $p\equiv29,101\pmod{120}$;\\
\indent \ \ \ \
$(1,2,3,2,1,3,2,3,1,2,1,3)$,  \ \ \ \ \ \ \ \ \ \ \ \ \ \ \ \ \ \ \ \ \ \ \ \ for $p\equiv53,77\pmod{120}$;\\
\indent \ \ \ \
$(1,2,4,2,4,6,2,-2,3,6,8,-2)$, \  \ \ \ \ \ \ \ \ \ \ \ \ \ \ \ \ \ \ for $p\equiv13,37\pmod{120}$;\\
\indent \ \ \ \
$(2,1,3,2,4,6,2,6,4,4,6,-2)$, \ \ \ \ \ \ \ \ \ \ \ \ \ \ \ \ \ \ \ \ \  for $p\equiv61,109\pmod{120}$.

\indent \ \ \ \   Let
$\B=\{(1,1,c_i)\cdot B_i\}$, where $c_i\in C_0^2/\{-1,1\}$ for $i=1,2$,
$c_i\in Z_p\setminus\{0\}$ for $i=3,4$, and $c_i\in C_0^2$ for
$5\leq i\leq 11$, then $\B$ forms a
$(4\times24p,4\times24,\{4,5\},$ $1)$-BDP by Lemma \ref{le33}. \qed

\section{Recursive Constructions}

\indent \ \ \ \ In this section, we shall give some recursive
constructions for $(u\times v,K,1)$-BDPs. We first introduce
some auxiliary designs.

\indent \ \ \ \ Let $G$ be an additive group of order $v$. A
$(G,k;1)$ difference matrix ($(G,k;1)$-DM for short) is a $k\times
v$ matrix $D=(d_{ij})$, $0\leq i\leq k-1$, $0\leq j\leq v-1$, with
entries from $G$, such that for any $0\leq i<j\leq k-1$, the
multiset $[d_{il}-d_{jl}| 0\leq l\leq v-1]$ contains every element
of $G$ exactly once. We usually write a $(G,k;1)$ difference matrix
as $(v,k;1)$-CDM if $G=Z_v$.



\begin{lemma}\label{dm-5}{\rm (\cite{Bur2})}
If there exist both a $(G,k;1)$-DM and a $(G',k;1)$-DM, then there exists a
$(G\times G',k;1)$-DM.
\end{lemma}

\begin{lemma}\label{dm-2}{\rm (\cite{Colbourn,GGN})}
Let $m$ and $k$ be positive integers such that {\rm
gcd}$(m,(k-1)!)=1$, then there exists an $(m,k;1)$-CDM. There exists
a $(3^a,5;1)$-CDM for any positive integer $a\geq3$.
\end{lemma}

\begin{lemma}\label{dm-4}{\rm (\cite{LZQW})}
If $a$, $b$ are positive integers, then there exists a
$(Z_{3^a}\times Z_{3^b},5;1)$-DM.
\end{lemma}

\begin{construction}{\rm (\cite{ZQ})}
\label{constru1} Suppose that both an $(u\times v,g\times
h,K,1)$-BDP and an optimal  $(g\times h,K,1)$-BDP exist, then
an optimal $(u\times v,K,1)$-BDP exists. Moreover, if the
given $(g\times h,K,1)$-BDP is a  $(g\times
h,s\times t,K,1)$-BDP, then so is the derived  $(u\times
v,K,1)$-BDP.
\end{construction}

\indent \ \ \ \ For a set $K$ of positive integers, let $k_{{\rm max}}=$max$\{k| k\in K\}$.

\begin{construction}{\rm (\cite{ZQ})}
\label{constru5} Suppose that there exist: ${\rm (1)}$ an
$(u\times v,g\times h,K,1)$-BDP; ${\rm (2)}$ a $(Z_m\times
Z_n,k_{{\rm max}};1)$-DM; ${\rm (3)}$ an $(mg\times nh,s\times
t,K,1)$-BDP. Then there exist an $(mu\times nv,s\times
t,K,1)$-BDP and an $(mu\times nv,mg\times nh,$ $K,1)$-BDP.
\end{construction}

\indent \ \ \ \  From the above constructions, it is not difficult
to obtain the following result.

\begin{corollary}
\label{constru2} Suppose that there exist: ${\rm (1)}$ an
$(u\times v,g\times h,K,1)$-BDP; ${\rm (2)}$ an $(m,k_{{\rm max}};1)$-CDM;
${\rm (3)}$ a $(g\times mh,s\times t,K,1)$-BDP (or
$(mg\times h,s\times t,K,1)$-BDP). Then there exist an
$(u\times mv,s\times t,K,1)$-BDP (or $(mu\times v,s\times
t,K,1)$-BDP) and an $(u\times mv,g\times mh,$ $K,1)$-BDP
(or $(mu\times v,mg\times h,$ $K,1)$-BDP).
\end{corollary}

\begin{corollary}\label{tl2}
If there exists an $(u\times v,g\times h,\{4,5\},1)$-BDP,
then there exists a  $(2u\times 6v, 2g\times
6h,\{4,5\},1)$-BDP.
\end{corollary}
{\bf Proof} \  A $(Z_2\times Z_6,5;1)$-DM is from Theorem 3.44 of
\cite{ACD}, the conclusion comes from Construction \ref{constru5}.
\qed

\section{The Proof of Theorem \ref{main1}}

In this section, we will prove Theorem \ref{main1}.

\begin{lemma}\label{lm51}
There exists an $(u\times v,g\times h,\{4,5\},1)$-BDP for
$(u,v,g,h)=(2,36,2,4), (2,72,2,8),$\\
$(2,108,2,12),(4,72,4,8),(6,12,2,4),(6,36,2,12),(12,$
$24,4,8),(18,4,2,4),$ $(18,12,2,12)$.
\end{lemma}
{\bf Proof} \ The blocks of the  desired  $(u\times v,g\times
h,\{4,5\},1)$-BDPs are displayed in Appendix A. \qed

\indent \ \ \ \  In this section, for positive integer $a$, define
$f(a)=\left \{
\begin{array}{ll}
8, & {\rm if}\  a \ \rm{is} \ even;\\
24, & {\rm if} \  a \ \rm{is} \ odd.
\end{array}
\right.$

\begin{lemma}\label{lm52}
If $a\geq0$ is an integer, then there exists a
$(4\times8\cdot3^a,4\times f(a),\{4,5\},1)$-BDP.
\end{lemma}
{\bf Proof} \ {\it 1st case: $a=0,1$}. The conclusion is trivial.

\indent \ \ \ \ {\it 2nd case: $a=2$}.  The conclusion is from Lemma
\ref{lm51}.

\indent \ \ \ \ {\it 3rd case: $a=3,4$}.

\indent \ \ \ \ A $(2\times4\cdot3^{a-1},2\times
4\cdot3^{a-3},\{4,5\},1)$-BDP is from Lemma \ref{lm51}, by applying
Corollary \ref{tl2}, we have a $(4\times8\cdot3^{a},4\times
8\cdot3^{a-2},\{4,5\},1)$-BDP. When $a=3$, the conclusion is
obtained. When $a=4$, the conclusion comes from Construction
\ref{constru1} by using a $(4\times72,4\times 8,\{4,5\},1)$-BDP.

\indent \ \ \ \ {\it 4th case: $a\geq5$}.

\indent \ \ \ \ Write $a=3n+a_1$, $n\geq1$ and $a_1\in\{2,3,4\}$.
Use induction on $n$. When $n=1$, a
$(4\times8\cdot3^{a_1},$ $4\times f(a_1),\{4,5\},1)$-BDP exists from the
above and a $(27,5;1)$-CDM exists from Lemma \ref{dm-2}, the
conclusion comes from Corollary \ref{constru2} by using a
$(4\times27f(a_1),4\times f(3+a_1),\{4,5\},1)$-BDP. Suppose that
there exists a $(4\times8\cdot3^{3(n-1)+a_1},4\times
f(3(n-1)+a_1),\{4,5\},1)$-BDP, then the conclusion is obtained from
Corollary \ref{constru2} by using a
$(4\times27f(3(n-1)+a_1),4\times f(3n+a_1),\{4,5\},1)$-BDP and a
$(27,5;1)$-CDM.\qed

\begin{lemma}\label{lm53}
If $a$ is a positive integer, then there exists a
$(12\times8\cdot3^a,4\times f(a+1),\{4,5\},1)$-BDP.
\end{lemma}
{\bf Proof} \ {\it 1st case: $a=1$}. The conclusion comes from Lemma
\ref{lm51}.

\indent \ \ \ \ {\it 2nd case: $a=2$}.

\indent \ \ \ \ A   $(6\times12,2\times 4,\{4,5\},1)$-BDP is
from Lemma \ref{lm51}, the conclusion is obtained from Corollary
\ref{tl2}.

\indent \ \ \ \ {\it 3rd case: $a=3$}.

\indent \ \ \ \ A  $(6\times36,2\times 12,\{4,5\},1)$-BDP is
from Lemma \ref{lm51}, by applying Corollary \ref{tl2}, we have a
 $(12\times216,$ $4\times 72,\{4,5\},1)$-BDP. So, the conclusion
comes from Construction \ref{constru1} by using a  $(4\times
72,4\times 8,\{4,5\},1)$-BDP in Lemma \ref{lm51}.

\indent \ \ \ \ {\it 4th case: $a\geq4$}.

\indent \ \ \ \ A  $(4\times8\cdot3^{a-2},4\times
f(a-2),\{4,5\},1)$-BDP is from Lemma \ref{lm52}, and a $(Z_3\times
Z_9,5;1)$-DM is from Lemma \ref{dm-4}, by applying Construction
\ref{constru5} with $u=4$, $v=8\cdot3^{a-2}$, we have a
$(12\times8\cdot3^a,12\times 9f(a-2),\{4,5\},1)$-BDP. If $a$ is odd,
then $f(a-2)=24$, so  the conclusion comes from Construction
\ref{constru1} by using a  $(12\times8\cdot 3^3,4\times
8,\{4,5\},1)$-BDP from the above. If $a$ is even, then $f(a-2)=8$, so
the conclusion comes from Construction \ref{constru1} by using a
 $(12\times8\cdot 3^2,4\times 24,\{4,5\},1)$-BDP from the
above. \qed

\begin{lemma}\label{lm54}
If $a$ is a positive integer, then there exists a
$(36\times8\cdot3^a,4\times f(a),\{4,5\},1)$-BDP.
\end{lemma}
{\bf Proof} \ {\it 1st case: $a=1$}.

\indent \ \ \ \ A   $(18\times4,2\times 4,\{4,5\},1)$-BDP is
from Lemma \ref{lm51}, the conclusion is obtained from Corollary
\ref{tl2}.

\indent \ \ \ \ {\it 2nd case: $a=2$}.

\indent \ \ \ \ A   $(18\times12,2\times 12,\{4,5\},1)$-BDP is
from Lemma \ref{lm51}, we have a   $(36\times72,4\times
72,\{4,5\},1)$-BDP from Corollary \ref{tl2},  the conclusion is
obtained from Construction \ref{constru1} by using a
$(4\times72,4\times 8,\{4,5\},1)$-BDP in Lemma \ref{lm51}.

\indent \ \ \ \ {\it 3rd case: $a\geq3$}.

\indent \ \ \ \ A $(4\times8\cdot3^{a-1},4\times
f(a-1),\{4,5\},1)$-BDP is from Lemma \ref{lm52}, and a $(Z_9\times
Z_3,5;1)$-DM is from Lemma \ref{dm-4}. Applying Construction
\ref{constru5} with $u=4$ and $v=8\cdot3^{a-1}$, we have a
$(36\times8\cdot3^a,36\times 3f(a-1),\{4,5\},1)$-BDP. If $a$ is odd,
then $f(a-1)=8$, so  the conclusion comes from Construction
\ref{constru1} by using a $(36\times24,4\times 24,\{4,5\},1)$-BDP
from the above. If $a$ is even, then $f(a-1)=24$, so the conclusion
comes from Construction \ref{constru1} by using a
$(36\times72,4\times 8,\{4,5\},1)$-BDP from the above. \qed

\begin{lemma}\label{lm55}
If $a, b\geq0$ are integers and $m,n$ are positive integers, then
there exists a   $(4\cdot3^a\times8\cdot3^b, $ $g\times
h,\{4,5\},1)$-BDP, where

$(g,h)=\left \{
\begin{array}{ll}
(4,8), & (a,b)=(0,0), (0,2n), (2m,0), (2m,2n), (2m-1,2n-1);\\
(4,24), &  (a,b)=(0,2n-1),(2m,2n-1),(2m-1,2n);\\
(12,8), & (a,b)=(2m-1,0).
\end{array}
\right.$
\end{lemma}
{\bf Proof} \ {\it 1st case: $a=0$ and $b\geq0$}.  The conclusion
comes from Lemma \ref{lm52}.

\indent \ \ \ \ {\it 2nd case: $a\geq1$ and $b=0$}.

\indent \ \ \ \ A   $(4\cdot3^a\times8,\frac{f(a)}{2}\times
8,\{4,5\},1)$-BDP which is equivalent to a
$(4\times8\cdot3^a,4\times f(a),\{4,5\},1)$-BDP is from Lemma
\ref{lm52}.

\indent \ \ \ \ {\it 3rd case: $a,b\geq1$}.

\indent \ \ \ \ A $(4\cdot3^{a-1}\times8,\frac{f(a-1)}{2}\times
8,\{4,5\},1)$-BDP which is equivalent to a
$(4\times8\cdot3^{a-1},4\times f(a-1),\{4,5\},1)$-BDP is from Lemma
\ref{lm52}. A $(Z_3\times Z_{3^b},5;1)$-DM is from Lemma \ref{dm-4},
we have a
$(4\cdot3^a\times8\cdot3^b,\frac{3f(a-1)}{2}\times
8\cdot3^b,\{4,5\},1)$-BDP  from Construction \ref{constru5}. If
$a$ is odd, we have $f(a-1)=8$, then the conclusion is obtained from
Construction \ref{constru1} by using a $(12\times8\cdot3^b,4\times
f(b+1),\{4,5\},1)$-BDP in Lemma \ref{lm53}. If $a$ is even, we have
$f(a-1)=24$, then the conclusion is obtained from Construction
\ref{constru1} by using a $(36\times8\cdot3^b,4\times
f(b),\{4,5\},1)$-BDP in Lemma \ref{lm54}. \qed

\begin{lemma}\label{lm55-1}
If $a, b\geq0$ are integers, then there exists a
$(4\cdot5^a\times24\cdot5^b,g\times h,\{4,5\},1)$-BDP, where
$(g,h)=\left \{
\begin{array}{ll}
(4,24), & a=b=0;\\
(4,120), &  a=0, b\geq1;\\
(20,24), & a\geq1, b\geq0.
\end{array}
\right.$

\end{lemma}
{\bf Proof} \ {\it 1st case: $a=b=0$}. The conclusion is
trivial.

\indent \ \ \ \ {\it 2nd case: $a=0$ and $b\geq1$}.

\indent \ \ \ \  A  $(4\times 40,4\times
8,\{4,5\},1)$-BDP exists from Lemma \ref{lm25}, and a $(15,5;1)$-CDM exists from
Example 2.1 of \cite{MC}, thus we have a  $(4\times
24\cdot5^2,4\times 120,\{4,5\},1)$-BDP from Corollary \ref{constru2}.
 Applying Corollary \ref{constru2} with a  $(4\times 24\cdot5^2,4\times
120,\{4,5\},1)$-BDP and a
$(5,5;1)$-CDM, one can obtain a  $(4\times 24\cdot5^3,4\times
120,\{4,5\},1)$-BDP. The  $(4\times
24\cdot5^b,4\times 120,\{4,5\},1)$-BDPs can be obtained by
repeating the process.

\indent \ \ \ \ {\it 3rd case: $a\geq1$ and $ b\geq0$}.

\indent \ \ \ \ If $a\geq1$ and $b=0$. A
$(4\cdot5^a\times24,20\times 24,\{4,5\},1)$-BDP which is equivalent
to a  $(4\times24\cdot5^a,4\times 120,\{4,5\},1)$-BDP is from
the above. If $a\geq1$ and $b\geq1$, a $(Z_5\times Z_3,5;1)$-DM
which is equivalent to a $(15,5;1)$-CDM exists from Example 2.1 of
\cite{MC}, we have a  $(20\times 120,20\times
24,\{4,5\},1)$-BDP from Construction \ref{constru5} by using a
 $(4\times 40,4\times 8,\{4,5\},1)$-BDP. There exist a
 $(4\times 24\cdot5^b,4\times 120,\{4,5\},1)$-BDP and a
$(5,5;1)$-CDM, one can obtain a  $(20\times 24\cdot5^b,20\times
24,\{4,5\},1)$-BDP from Corollary \ref{constru2} by using a
$(20\times 120,20\times 24,\{4,5\},1)$-BDP. Since there exist a
$(5^b,5;1)$-CDM from Lemma \ref{dm-2} and a
$(4\cdot5^a\times 24,20\times 24,\{4,5\},1)$-BDP, then a
$(4\cdot5^a\times 24\cdot5^b,20\times 24,\{4,5\},1)$-BDP  exists from
Corollary \ref{constru2} by using a  $(20\times
24\cdot5^b,20\times 24,\{4,5\},1)$-BDP.  \qed

\begin{lemma}\label{lm56}
If $u$ and $v$ are positive integers such that {\rm gcd}$(uv,6)=1$,
then there exists a  $(4u\times 8v, 4\times 8,$ $\{4,5\},1)$-BDP.
\end{lemma}
{\bf Proof} \ {\it 1st case: $u=v=1$}. The conclusion is
trivial.

\indent \ \ \ \ {\it 2nd case: $u=1$ and $v>1$}. Let
$v=p_1^{a_1}p_2^{a_2}\ldots p_l^{a_l}$ be the factorization of $v$,
where each $p_i\geq5$ be prime and each integer $a_i\geq1$. For each
prime $p_i$,   a $(p_i,5;1)$-CDM exists from Lemma \ref{dm-2}, and
a  $(4\times 8p_i,4\times 8,\{4,5\},1)$-BDP exists from Lemma
\ref{lm25}.   Applying Corollary \ref{constru2} with a  $(4\times 8p_1,4\times
8,\{4,5\},1)$-BDP and a $(p_i,5;1)$-CDM, one can obtain a  $(4\times 8p_1p_i,4\times
8p_i,\{4,5\},1)$-BDP. By using Construction \ref{constru1} with a
 $(4\times 8p_i,4\times 8,\{4,5\},1)$-BDP, one can obtain a
 $(4\times 8p_1p_i,4\times 8,\{4,5\},1)$-BDP. Repeat the
process, one can get a  $(4\times 8v,4\times 8,\{4,5\},1)$-BDP.

\indent \ \ \ \ {\it 3rd case: $u>1$ and $v=1$}. Since a
$(4u\times 8,4\times 8,\{4,5\},1)$-BDP is equivalent to a
$(4\times 8u,4\times 8,\{4,5\},1)$-BDP,  the conclusion is from
the 2nd case.

\indent \ \ \ \ {\it 4th case: $u>1$ and $v>1$}.

\indent \ \ \ \ A  $(4u\times 8,4\times
8,\{4,5\},1)$-BDP exists from the above, and a $(v,5;1)$-CDM is from Lemma
\ref{dm-2},  the conclusion comes from Corollary \ref{constru1} by
using a  $(4\times 8v,4\times 8,\{4,5\},1)$-BDP. \qed

\begin{lemma}\label{lm56-2}
If $u,v$ are positive integers such that {\rm gcd}$(uv,30)=1$, then
there exists a  $(gu\times hv,g\times h,$ $\{4,5\},1)$-BDP for
$(g,h)=(4,24),(4,120),(12,8),$ $(20,24)$.
\end{lemma}
{\bf Proof} \ Since a  $(gu\times hv,g\times
h,\{4,5\},1)$-BDP is equivalent to a  $(g'u\times
h'v,g'\times h',\{4,5\},1)$-BDP,
where $(g,h,g',h')=(4,24,12,8)$, or $(4,120,20,24)$,  then we need only to consider
the cases of $(g,h)=(4,24),(4,120)$.\\
\indent \ \ \ \ Let $v=p_1^{a_1}p_2^{a_2}\ldots p_l^{a_l}$ be the
factorization of $v$,  $p_i>5$ be prime and
$a_i\geq1$, $1\leq i\leq l$.

{\it 1st case: $(g,h)=(4,24)$}.

\indent \ \ \ \  For each  $p_i$,  a $(p_i,5;1)$-CDM exists from
Lemma \ref{dm-2}, and a  $(4\times 24p_i,4\times
24,\{4,5\},1)$-BDP exists from Lemmas \ref{lm26} and \ref{lm27}. Similar to
 the proof of Lemma \ref{lm56}, one can get the conclusion.

{\it 2nd case: $(g,h)=(4,120)$}.

\indent \ \ \ \  A  $(4u\times 24v,4\times
24,\{4,5\},1)$-BDP exists from the above,  the conclusion comes
 from Corollary \ref{constru2} by using a $(5,5;1)$-CDM in Lemma
\ref{dm-2}.  \qed

\indent \ \ \ \ Now we are in a position to prove Theorem
\ref{main1}.

{\bf The Proof of Theorem \ref{main1}} \ Write $u=3^au_1$,
$v=3^bv_1$, where {\rm gcd}$(u_1v_1,6)=1$, and $a,b\geq0$.
Then $4u\times8v=4\cdot3^au_1\times8\cdot3^bv_1$.

\indent \ \ \ \ {\it 1st case: $u_1=v_1=1$}.

\indent \ \ \ \  A  $(4\cdot3^a\times8\cdot3^b,g\times
h,\{4,5\},1)$-BDP exists from Lemma \ref{lm55}, where $(g,h)=(4,8)$,
$(4,24)$, or $(12,8)$.  An
optimal   $(4\times24,\{4,5\},1)$-BDP  is equivalent to
an optimal  $(12\times8,\{4,5\},1)$-BDP, the blocks of an optimal
  $(4\times24,\{4,5\},1)$-BDP is listed below.

\{(0,0),(0,1),(0,3),(0,7),(1,0)\},
\{(0,0),(0,5),(0,13),(1,1),(1,11)\},

\{(0,0),(0,9),(1,3),(2,5)\}, \{(0,0),(1,4),(2,9),(3,17)\}.

The conclusion comes from Construction
\ref{constru1}.

\indent \ \ \ \ {\it 2nd case: $u_1=1$ and $v_1>1$}.

\indent \ \ \ \  Let $v_1=5^cv_2$, where ${\rm gcd}(v_2,30)=1$. If $c=0$.
A  $(4\cdot3^a\times8\cdot3^b,g\times h,\{4,5\},1)$-BDP
exists from Lemma \ref{lm55}, where $(g,h)=(4,8)$, $(4,24)$, or
$(12,8)$, and a $(v_1,5;1)$-CDM exists from Lemma \ref{dm-2}, we
have a  $(4\cdot3^a\times8\cdot3^bv_1,g\times
hv_1,\{4,5\},1)$-BDP from Corollary \ref{constru2}. Since there exist
a  $(12\times8v_1,12\times 8,\{4,5\},1)$-BDP and a
$(4\times24v_1,4\times 24,\{4,5\},1)$-BDP from Lemma \ref{lm56-2}, a
 $(4\times8v_1,4\times 8,\{4,5\},1)$-BDP is from Lemma
\ref{lm56}, we have a  $(4\cdot3^a\times8\cdot3^bv_1,g\times
h,\{4,5\},1)$-BDP from Construction \ref{constru1}. So, the
conclusion is obtained from Construction \ref{constru1},
where the  optimal  $(g\times
h,\{4,5\},1)$-BDP, $(g,h)=(4,24),(12,8)$ is from the above.

\indent \ \ \ \ If $c\geq1$.  A
$(4\cdot3^a\times8\cdot3^bv_2,g\times h,\{4,5\},1)$-BDP exists from the
above,  where $(g,h)=(4,8)$, $(4,24)$, or $(12,8)$,  and a
$(5^c,5;1)$-CDM exists from Lemma \ref{dm-2}, we have a
$(4\cdot3^a\times8\cdot3^b5^cv_2,g\times 5^ch,\{4,5\},1)$-BDP. Since
a  $(4\times24\cdot 5^c,4\times 120,\{4,5\},1)$-BDP, which is
equivalent to a  $(12\times8\cdot 5^c,12\times
40,\{4,5\},1)$-BDP, exists from Lemma \ref{lm55-1}, and a
$(4\times8\cdot 5^c,4\times 8,\{4,5\},1)$-BDP is from Lemma
\ref{lm56}, then we obtain a
$(4\cdot3^a\times8\cdot3^b5^cv_2,g'\times h',\{4,5\},1)$-BDP, where
$(g',h')=$(4,8), (4,120) or (12,40). So, the conclusion is obtained
from Construction \ref{constru1}, where the
 optimal  $(4\times 120,\{4,5\},1)$-BDP
 (equivalent to an optimal  $(12\times 40,\{4,5\},1)$-BDP) is listed in
Appendix B.

\indent \ \ \ \ {\it 3rd case: $u_1>1$ and $v_1=1$}.

\indent \ \ \ \ An optimal   $(4\cdot3^au_1\times
8\cdot3^b,\{4,5\},1)$-BDP, which is equivalent to an optimal
$(4\cdot3^a\times 8\cdot3^bu_1,\{4,5\},1)$-BDP, exists from the second
case.

\indent \ \ \ \ {\it 4th case: $u_1>1$ and $v_1>1$}.

\indent \ \ \ \ Let $u_1=5^cu_2$, $v_1=5^dv_2$, where
${\rm gcd}(u_2v_2,30)=1$. If $c=d=0$. An $(u_2,5;1)$-CDM and a
$(v_2,5;1)$-CDM exist, we have a $(Z_{u_2}\times Z_{v_2},5;1)$-DM
from Lemma \ref{dm-5}. There exists a   $(gu_2\times
hv_2,g\times h,\{4,5\},1)$-BDP from Lemmas \ref{lm56}, \ref{lm56-2},
where $(g,h)=(4,8)$, $(4,24)$, or $(12,8)$, we have a
$(4\cdot3^au_2\times8\cdot3^bv_2,g\times h,\{4,5\},1)$-BDP from
Construction \ref{constru5} by using a
$(4\cdot3^a\times8\cdot3^b,g\times h,\{4,5\},1)$-BDP and a
$(Z_{u_2}\times Z_{v_2},5;1)$-DM, and hence the conclusion is from
Construction \ref{constru1}.

\indent \ \ \ \ If $c=0, \ d\geq1$. A $(5^d,5;1)$-CDM
exists from Lemma \ref{dm-2}, and a
$(4\cdot3^au_2\times8\cdot3^bv_2,g\times h,\{4,5\},1)$-BDP exists from the
above, similar to the proof of the second case, one can obtain
the conclusion.

\indent \ \ \ \  If $c\geq1, \ d=0$. This case is equivalent to the
 case of $c=0, \ d\geq1$.

\indent \ \ \ \  If $c\geq1, \ d\geq1$. A $(5^c,5;1)$-CDM and a
$(5^d,5;1)$-CDM exist, we have a $(Z_{5^c}\times Z_{5^d},5;1)$-DM
from Lemma \ref{dm-5}, then  a
$(4\cdot3^a5^cu_2\times8\cdot3^b5^dv_2,5^cg\times
5^dh,\{4,5\},1)$-BDP exists  from Construction \ref{constru5} by using a
  $(4\cdot3^au_2\times8\cdot3^bv_2,g\times h,\{4,5\},1)$-BDP,
where $(g,h)=(4,8)$, $(4,24)$, or $(12,8)$. Since there exist a
  $(4\cdot5^c\times8\cdot5^d,4\times 8,\{4,5\},1)$-BDP from
Lemma \ref{lm56}, and a   $(4\cdot5^c\times24\cdot5^d,20\times
24,\{4,5\},1)$-BDP (equivalent to a
$(12\cdot5^c\times8\cdot5^d,60\times 8,\{4,5\},1)$-BDP) from Lemma
\ref{lm55-1}, we get a
$(4\cdot3^a5^cu_2\times8\cdot3^b5^dv_2,g'\times h',\{4,5\},1)$-BDP
from Construction \ref{constru1}, where $(g',h')=$(4,8), (20,24), or
(60,8). So, the conclusion is obtained from Construction
\ref{constru1}, where the
optimal   $(20\times 24,\{4,5\},1)$-BDP and  optimal
  $(60\times 8,\{4,5\},1)$-BDP exist since they are  equivalent
to the optimal   $(4\times 120,\{4,5\},1)$-BDP in Appendix B. \qed

\section{Application to OOSPCs}

\indent \ \ \ \  Optical orthogonal signature pattern code (OOSPC) was introduced
by Kitayama \cite{K} for parallel transmission of 2-D image in
multicore fiber optical code division multiple access (OCDMA)
networks. Multiple-weight (MW) OOSPCs were introduced by Kwong and Yang
to meet multiple quality of services (QoS) requirement \cite{KY1}. For OOSPCs and MW-OOSPCs,
the interested readers may refer to \cite{CJL}, \cite{CJL1}, \cite{JDWG}, \cite{K},
\cite{LZQW},   \cite{PC}, \cite{PC1},
\cite{PC2}, \cite{PC3}, \cite{S}, \cite{SK},  \cite{YK},
\cite{ZQ} and the references therein.

\indent \ \ \ \  The following result is an analogue
result (Theorem II.1) in \cite{S}.

\begin{lemma}\label{th1.1}
An optimal balanced  $(u, v, K, 1)$-OOSPC is equivalent to an optimal
$(u\times v, K, 1)$-BDP.
\end{lemma}

\indent \ \ \ \ By using Lemma \ref{th1.1} and Theorem \ref{main1},
the following result is  obtained.

\begin{corollary}\label{CM}
There exist optimal balanced $(4u,8v, \{4,5\},1)$-OOSPCs for any
odd integers $u>1, \ v>1$.
\end{corollary}

\noindent {\bf Acknowledgements} \ The authors  would  like to thank
 Professor Gennian Ge  of Capital Normal University,  China, for
 helpful discussions. The first author is supported in part  by
Guangxi Nature Science Foundation (No. 2018GXNSFBA138038) and BAGUI
Scholar Program of Guangxi Zhuang Autonomous Region of China (No.
201979). The second author is supported by Guangxi Nature Science
Foundation (No. 2019GXNSFBA245021). The last author is supported in
part by NSFC (No. 11671103, 11801103) and Guangxi Nature Science
Foundation (No. 2017GXNSFBA198030). \\

 \noindent  {\bf Appendix A} \

\noindent  $(u,v,g,h)=${\bf (2,36,2,4)}

\{(0,0),(1,14),(1,15),(1,26)\}, \{(0,0),(0,5),(0,15),(0,19),(1,7)\},
\{(0,0),(0,2),(0,8),(1,13)\},

\{(0,0),(0,3),(1,4),(1,20),(1,33)\}.

\noindent $(u,v,g,h)=${\bf (2,72,2,8)}

\{(0,0),(0,71),(1,1),(1,14)\},
\{(0,0),(0,17),(0,37),(0,48),(1,20)\}, \{(0,0),(0,4),(0,7),(0,53)\},

\{(0,0),(0,33),(0,38),(1,22),(1,43)\},
 \{(0,0),(0,8),(0,30),(1,42)\},
\{(0,0),(0,14),(0,70),(1,46)\},

\{(0,0),(0,43),(1,4),(1,51),(1,66)\}, \{(0,0),(0,28),(0,60),(1,35),(1,41)\}.

\noindent $(u,v,g,h)=${\bf (2,108,2,12)}

\{(0,0),(0,65),(0,89),(1,52)\},
\{(0,0),(0,5),(0,37),(1,43),(1,49)\},
\{(0,0),(0,30),(0,47),(0,68)\},

\{(0,0),(0,59),(0,60),(0,93),(1,91)\},
\{(0,0),(0,10),(0,52),(0,83),(1,33)\},

\{(0,0),(0,46),(0,96),(1,35),(1,39)\},
\{(0,0),(0,3),(0,23),(1,28),(1,107)\},

\{(0,0),(0,13),(1,53),(1,92)\}, \{(0,0),(0,8),(0,22),(1,30)\},
\{(0,0),(0,51),(0,67),(1,93)\},

\{(0,0),(0,7),(1,10),(1,21),(1,74)\}, \{(0,0),(0,2),(0,28),(1,48)\}.

\noindent $(u,v,g,h)=${\bf (4,72,4,8)}

\{(0,0),(0,24),(1,16),(2,17),(2,64)\},
\{(0,0),(1,49),(2,16),(2,35)\}, \{(0,0),(1,11),(3,60),(3,64)\},

\{(0,0),(0,65),(0,70),(3,4),(3,42)\},
\{(0,0),(0,16),(2,6),(3,26),(3,37)\},

\{(0,0),(1,53),(1,70),(3,39),(3,40)\},
\{(0,0),(1,15),(1,21),(2,25),(3,67)\},

\{(0,0),(1,17),(2,60),(3,13)\},
\{(0,0),(0,43),(1,56),(1,69),(2,44)\},
\{(0,0),(0,49),(2,11),(2,51)\},

\{(0,0),(0,51),(1,19),(1,29),(3,58)\},
\{(0,0),(1,55),(3,5),(3,31)\}, \{(0,0),(0,58),(1,24),(1,57)\},

\{(0,0),(0,15),(0,37),(3,35),(3,65)\}, \{(0,0),(0,8),(0,20),(2,5)\},
\{(0,0),(0,28),(0,31),(1,34)\}.

\noindent $(u,v,g,h)=${\bf (6,12,2,4)}

\{(0,0),(0,5),(1,8),(3,4)\}, \{(0,0),(1,11),(2,4),(3,5),(5,6)\},
\{(0,0),(0,1),(2,3),(4,3)\},

\{(0,0),(0,2),(1,0),(1,4),(5,5)\}.

\noindent $(u,v,g,h)=${\bf (6,36,2,12)}

\{(0,0),(1,5),(3,20),(5,11)\}, \{(0,0),(3,26),(4,24),(4,25)\},
\{(0,0),(0,16),(3,2),(4,8),(4,16)\},

\{(0,0),(1,17),(3,7),(4,34)\}, \{(0,0),(2,5),(2,19),(3,5)\},
\{(0,0),(1,29),(2,31),(3,25),(3,35)\},

\{(0,0),(0,32),(1,20),(2,33),(2,35)\},
\{(0,0),(0,7),(4,27),(5,24)\}, \{(0,0),(1,11),(2,7),(3,28)\},

\{(0,0),(1,16),(2,23),(2,34),(4,12)\}, \{(0,0),(1,9),(1,28),(3,13),(5,35)\},

\{(0,0),(0,5),(1,8),(1,31),(2,18)\}.

\noindent $(u,v,g,h)=${\bf (12,24,4,8)}

\{(0,0),(2,1),(2,20),(6,4),(10,21)\}, \{(0,0),(4,15),(6,5),(10,1)\},
\{(0,0),(2,22),(3,4),(7,20)\},

\{(0,0),(2,0),(5,10),(9,22),(9,23)\},
\{(0,0),(0,2),(1,3),(2,8),(11,10)\},

\{(0,0),(4,22),(5,9),(10,15),(11,1)\},
\{(0,0),(0,7),(1,22),(7,0),(8,17)\},

\{(0,0),(2,12),(7,7),(11,16)\}, \{(0,0),(1,0),(8,13),(11,20)\},
\{(0,0),(0,13),(1,20),(10,6)\},

\{(0,0),(2,19),(5,14),(9,16),(10,13)\},
\{(0,0),(1,2),(1,18),(3,23),(5,12)\},

 \{(0,0),(4,0),(5,13),(9,8)\},
\{(0,0),(2,10),(5,15),(7,6)\}, \{(0,0),(0,4),(3,17),(4,5)\},

\{(0,0),(5,16),(6,1),(6,11),(8,3)\}.

\noindent $(u,v,g,h)=${\bf (18,4,2,4)}

\{(0,0),(1,1),(8,3),(11,3),(14,1)\}, \{(0,0),(1,0),(6,3),(11,0)\},
\{(0,0),(3,3),(6,0),(7,3),(8,1)\},

 \{(0,0),(2,0),(4,2),(6,1)\}.

\noindent $(u,v,g,h)=${\bf (18,12,2,12)}

 \{(0,0),(4,11),(5,4),(15,6),(17,3)\},
\{(0,0),(7,7),(10,10),(14,3),(15,10)\},

\{(0,0),(1,11),(2,1),(13,9)\}, \{(0,0),(3,7),(11,3),(17,9)\},
\{(0,0),(8,1),(10,3),(11,11),(12,9)\},

\{(0,0),(2,10),(13,11),(15,3),(17,6)\},
\{(0,0),(3,1),(5,9),(15,4)\}, \{(0,0),(4,2),(10,0),(11,0)\},

\{(0,0),(2,0),(5,11),(8,4),(13,10)\}, \{(0,0),(1,4),(4,4),(5,5)\},
\{(0,0),(2,7),(6,7),(8,6),(12,0)\},

 \{(0,0),(2,5),(6,8),(13,2)\}.

\noindent  {\bf Appendix B} \
An optimal  $(4\times 120,\{4,5\},1)$-BDP.

\{(0,0),(1,103),(3,27),(3,73)\}, \{(0,0),(1,92),(1,98),(2,107)\},
\{(0,0),(0,4),(1,17),(2,35)\},

\{(0,0),(1,51),(2,45),(3,32)\}, \{(0,0),(1,46),(2,16),(3,113)\},
\{(0,0),(1,45),(1,76),(2,27)\},

\{(0,0),(0,15),(0,48),(2,88),(3,4)\},
\{(0,0),(0,88),(2,99),(3,61),(3,98)\},

\{(0,0),(0,69),(2,86),(3,50),(3,112)\},
\{(0,0),(0,25),(0,84),(3,34),(3,78)\},

\{(0,0),(0,113),(1,23),(1,101),(3,47)\},
\{(0,0),(1,25),(2,39),(3,51)\},

\{(0,0),(0,28),(0,57),(0,77),(0,98)\},
\{(0,0),(0,30),(0,54),(0,94),(3,33)\},

\{(0,0),(1,106),(3,58),(3,92)\}, \{(0,0),(0,9),(1,52),(2,4)\},
\{(0,0),(0,47),(2,98),(2,110)\},

\{(0,0),(1,53),(1,54),(2,55)\}, \{(0,0),(1,35),(2,91),(3,29)\},
\{(0,0),(1,10),(1,75),(2,70),(3,7)\},

\{(0,0),(0,23),(1,39),(1,78),(2,43)\},
\{(0,0),(0,8),(1,32),(2,95),(3,24)\},

\{(0,0),(0,3),(0,13),(2,41),(2,59)\},
\{(0,0),(0,17),(0,19),(1,48),(1,100)\},

\{(0,0),(0,14),(2,111),(3,55),(3,100)\},
\{(0,0),(0,11),(0,115),(2,113),(3,26)\},

\{(0,0),(1,80),(3,79),(3,117)\}, \{(0,0),(1,37),(2,12),(3,52)\}.

\end{document}